\documentclass[12pt]{amsart}

\usepackage{amsthm,amssymb,amscd}
\usepackage{epic,eepic}
\usepackage[all]{xy}
\usepackage[dvips]{graphicx}

\newcommand     {\comment}[1]   {}
\newcommand{\mute}[2] {}
\newcommand     {\printname}[1] {}

\newcommand{\labell}[1] {\label{#1}\printname{#1}}

\swapnumbers

\numberwithin{equation}{section}
\newtheorem {Theorem}                   {Theorem}
\newtheorem {Lemma}[Theorem]         	{Lemma}

\newtheorem{prop}[Theorem]              {Proposition}
\newtheorem{Proposition}[Theorem]       {Proposition}

\newtheorem*{Proposition*}              {Proposition}

\theoremstyle{definition}
\newtheorem{Definition}[Theorem]{Definition}
\newtheorem{noTitle}[Theorem]{}

\theoremstyle{remark}
\newtheorem{Remark}[Theorem]{Remark}
\newtheorem*{remark}{Remark}
\newtheorem*{warning}{Warning}
\newtheorem{Example}[Theorem]{Example}
\newtheorem*{Example*}{Example}

\def \: {\colon}

\def \calT {{\mathcal T}} 
\def \h {{\mathfrak h}} 
\def \t {{\mathfrak t}} 

\def \R {{\mathbb R}}
\def \C {{\mathbb C}}

\def \eps {\varepsilon}
\def \vareps {\varepsilon}
\def \half {{\frac{1}{2}}}
\def \ssminus {\smallsetminus}

\def \gammaRn {\overline{\gamma}}

\def    \inv    {^{-1}}

\def    \ol     {\overline}
\newcommand{\xhalf}{x_{1/2}}

\def \calT {{\mathcal T}}

\newcommand{\Uxeps}{U_{[x],\varepsilon}}
\newcommand{\Uyeps}{U_{[y],\varepsilon}}
\newcommand{\Wxeps}{W_{[x],\varepsilon}}
\def \tU {\tilde{U}}
\newcommand{\tUxeps}{\tU_{[x],\varepsilon}}

\begin{document}

\title[Local to Global Convexity] {Revisiting Tietze-Nakajima -- \\
local and global convexity for maps}

\author[Yael Karshon]{Yael Karshon}
\email{karshon@math.toronto.edu}

\author[Christina Bjorndahl]{Christina Bjorndahl ${}^*$}
\thanks{$^*$ formerly Christina Marshall}
\email{cjm295@cornell.edu}
\address{Department of Mathematics, University of Toronto, Toronto, ON 
M5S 2E4, CANADA}

\maketitle


\section{Introduction}

A theorem of Tietze and Nakajima, from 1928, 
asserts that if a subset $X$ of $\R^n$ is closed, connected, 
and locally convex, then it is convex \cite{T,N}.
There are many generalizations of this ``local to global convexity"
phenomenon in the literature; a partial list is 
\cite{BF,cel,kay,KW,klee,SSV,S,tam}.

This paper contains an analogous ``local to global convexity" theorem
when the inclusion map of $X$ to $\R^n$ is replaced by a map
from a topological space $X$ to $\R^n$ that satisfies
certain local properties:
We define a map $\Psi \colon X \to \R^n$ to be convex if 
any two points in $X$ can be connected by a path $\gamma$
whose composition with $\Psi$ parametrizes a straight line segment
in $\R^n$ and this parametrization is monotone along the segment.
%
%
See Definition \ref{def:convex}.
We show that, if $X$ is connected and Hausdorff, $\Psi$ is proper,
and each point has a neighbourhood $U$ such that $\Psi|_U$ is convex 
and open as a map to its image, then $\Psi$ is convex and open as
a map to its image. We deduce that the image of $\Psi$ is convex
and the level sets of $\Psi$ are connected.
See Theorems \ref{maintheorem} and \ref{Theorem}. 

Our motivation comes from the Condevaux-Dazord-Molino proof 
\cite{CDM,HNP} of the Atiyah-Guillemin-Sternberg convexity 
theorem in symplectic geometry \cite{A,GS}.  
See section~\ref{sec:moment}.

This paper is the result of an undergraduate research project
that spanned over the years 2004--2006.  
The senior author takes the blame for the delay in publication
after posting our arXiv eprint.  
While preparing this paper we learned of the paper \cite{BOR1}
by Birtea, Ortega, and Ratiu, which achieves similar goals.  
In section \ref{sec:ratiu} we discuss relationships between
our results and theirs.
After \cite{BOR1}, our results are not essentially new, 
but our notion of ``convex map" gives elegant statements,
and our proofs are so elementary that they are accessible 
to undergraduate students with basic topology background.

\subsection*{Acknowledgements}
The first author is partially supported by an NSERC
Discovery grant. 
The second author was partially funded by an NSERC USRA grant
in the summers of 2004 and 2005.
The authors are grateful to River Chiang and to Tudor Ratiu 
for helpful comments on the manuscript.

\section{The Tietze-Nakajima Theorem}

Let $B(x,r)$ [$\overline{B}(x,r)$]
denote the open [closed] ball in $\R^n$ of radius $r$, centered at $x$.  
A closed subset $X$ of $\R^n$ is \emph{locally convex} 
if for every $x \in X$ there exists $\delta_x > 0$ such that
$B(x, \delta_x) \cap X$ is convex.
The Tietze-Nakajima theorem \cite{T,N} asserts that ``local convexity
implies global convexity":

\begin{Theorem}[Tietze-Nakajima]
Let $X$ be a closed, connected, and locally convex subset of $\R^n$.
Then $X$ is convex.
\end{Theorem}

\begin{Example}
A disjoint union of two closed balls is closed and locally convex
but is not connected.  A punctured disk is connected
and satisfies the locally convexity condition but it is not closed.
\end{Example}

A closed subset $X \subset \R^n$ is \emph{uniformly locally convex} 
on a subset $A \subset X$ if there exists $\delta > 0$ such that 
$B(x, \delta) \cap X$ is convex for all $x \in A$.
%
%
\begin{Lemma} [Uniform local convexity on compact sets]
\labell{Lemma 1} 
Let $X$ be a closed subset of $\R^n$.  If $X$ is locally convex
and $A \subset X$ is compact, then $X$ is uniformly locally convex on $A$.
\end{Lemma}

\begin{proof}
Since $X$ is locally convex, for every $x \in X$ there exists a
$\delta_x > 0$ such that $B(x, \delta_x) \cap X$ is convex.  By
compactness there exist points $x_1, \ldots, x_n$ such that 
$A \subset \bigcup_{i=1}^{n} B(x_i,\frac{1}{2}\delta_{x_i})$.  
Let $\delta = \min\{\frac{1}{2}\delta_{x_i}\}$.  Then for every $x \in A$
there exists $i$ such that $B(x,\delta) \subset B(x_i,\delta_{x_i})$.  
It follows that $B(x,\delta) \cap X$ is convex.
\end{proof}

\begin{noTitle}
Let $X$ be a closed, connected, locally convex subset of $\R^n$.
For two points $x_0$ and $x_1$ in $X$, define their
\emph{distance in $X$}, denoted $d_X(x_0,x_1)$, as follows:
\begin{displaymath}
d_X(x_0,x_1) = \inf\{l(\gamma) \,|\, \gamma \colon \left[0,1\right]
\rightarrow X, \ \gamma(0)=x_0, \ \gamma(1)=x_1 \},
\end{displaymath}
where $l(\gamma)$ is the length of the path $\gamma$.
\end{noTitle}

In this definition it doesn't matter if we take the infimum 
over continuous paths or polygonal paths:
let $\gamma \colon [0,1] \to X$ be a continuous path in $X$.
Let $\delta$ be the radius associated with uniform local convexity on
the compact set $\{\gamma(t), 0 \leq t \leq 1 \}$.
By uniform continuity of $\gamma$ on the compact interval $[0,1]$,
there exist $0 = t_0 < t_1 < \ldots < t_k = 1$ such that 
$\|\gamma(t_{i-1}) - \gamma(t_i)\| < \delta$ for $i=1, \ldots, k$.  
The polygonal path through the points
$\gamma(t_0), \ldots, \gamma(t_{k})$ is contained in $X$
and has length~$\leq l(\gamma)$.

Also note that $d_X(x_0,x_1) \geq \| x_1 - x_0 \|$,
with equality if and only if the segment $[x_0,x_1]$ is contained in~$X$.

%
%

\begin{Lemma}[Existence of midpoint]
\labell{Lemma 3} 
Let $X$ be a closed, connected, and locally convex
subset of $\R^n$.  Let $x_0$ and $x_1$ be in $X$.  
Then there exists a point $x_{1/2}$ in $X$ such that 
\begin{equation}
\labell{midpoint}
d_X(x_0, x_{1/2}) = d_X(x_{1/2}, x_1) = \frac{1}{2}d_X(x_0,x_1).
\end{equation}
\end{Lemma}

\begin{proof}
Let $\gamma_j$ be paths in $X$ connecting $x_0$ and $x_1$ 
such that $\{l(\gamma_j)\}$ converges to $d_X(x_0, x_1)$.
Let $t_j \in [0,1]$ be such that $\gamma_j(t_j)$ is the midpoint of
the path $\gamma_j$:
\begin{equation*}
l(\gamma_j \arrowvert_{\scriptscriptstyle{\left[0,t_j\right]}}) = 
l(\gamma_j \arrowvert_{\scriptscriptstyle{\left[t_j,1\right]}}) =
\frac{1}{2}l(\gamma_j).
\end{equation*} 
Since the sequence of midpoints $\{\gamma_j(t_j)\}$ is bounded and $X$
is closed, this sequence has an accumulation point $\xhalf \in X$.  
We will show that the point $\xhalf$ satisfies equation~\eqref{midpoint}.

We first show that for every $\vareps > 0$ there exists a path $\gamma$
connecting $x_0$ and $x_\half$ such that 
$l(\gamma) < \half d_X(x_0,x_1) + \vareps$.

Let $\delta > 0$ be such that $B(x_\half,\delta) \cap X$ is convex.
Let $j$ be such that 
$ \| \gamma_j(t_j) - x_\half \| < \min ( \delta, \frac{\vareps}{2}) $
and such that $ l(\gamma_j) < d_X(x_0,x_1) + \vareps $.
The segment $[\gamma_j(t_j),x_\half]$ is contained in $X$.
Let $\gamma$ be the concatenation of $\gamma_j|_{[0,t_j]}$
with this segment.  Then $\gamma$ is a path in $X$
that connects $x_0$ and $x_\half$, 
and $l(\gamma) < \half d_X(x_0,x_1) + \vareps$.

Thus, $d_X(x_0,x_\half) \leq \half d_X (x_0,x_1)$.
By the same argument, $d_X(x_\half,x_1) \leq \half d_X(x_0,x_1)$.
If either of these were a strict inequality, then it would be
possible to construct a path in $X$ from $x_0$ to $x_1$ whose length 
is less than $d_X(x_0,x_1)$, which contradicts the definition 
of $d_X(x_0,x_1)$.
\end{proof}

%
%
\begin{proof}[Proof of the Tietze-Nakajima theorem]
Fix $x_0$ and $x_1$ in $X$.

By Lemma \ref{Lemma 3}, there exists a point $x_{1/2}$ such that
\begin{displaymath}
d_X(x_0, x_{1/2}) = d_X(x_{1/2}, x_1) = \frac{1}{2}d_X(x_0, x_1).
\end{displaymath}
Likewise, there exists a point $x_{1/4}$ that satisfies
\begin{displaymath}
d_X(x_0, x_{1/4}) = d_X(x_{1/4}, x_{1/2}) = \frac{1}{2}d_X(x_0,
x_{1/2}).
\end{displaymath}
By iteration, we get a map $\frac{j}{2^m} \mapsto x_{\frac{j}{2^m}}$,
for nonnegative integers $j$ and $m$ where $0 \leq j \leq 2^m$, 
such that
\begin{equation}
\labell{iteratedmidpoint}
d_X(x_\frac{j-1}{2^m}, x_\frac{j}{2^m}) = d_X(x_\frac{j}{2^m},
x_\frac{j+1}{2^m}) = \frac{1}{2}d_X(x_\frac{j-1}{2^m},
x_\frac{j+1}{2^m}).
\end{equation}

Let
$$ r > d_X(x_0,x_1).$$ 
For all $0 \leq j \leq 2^m$, the following is true:
\begin{equation*}
\|x_{\frac{j}{2^m}} - x_0\| \leq d_X(x_{\frac{j}{2^m}}, x_0) \leq
\sum_{i=1}^j d_X(x_{\frac{i-1}{2^m}}, x_{\frac{i}{2^m}}) =
\frac{j}{2^m}d_X(x_0, x_1) < r.
\end{equation*}
Thus $x_{\frac{j}{2^m}}$ belongs to the compact set
\begin{equation*}
\overline{B}(x_0,r)\cap X.
\end{equation*}
Let $\delta$ denote the radius associated with uniform local convexity 
on this compact set.  Choose $m$ large enough such that 
$\frac{1}{2^m}d_X(x_0, x_1) < \delta$.  
Since the intersection
$B(x_{\frac{j}{2^m}}, \delta) \cap X$ is convex
and $x_{\frac{j-1}{2^m}} \in B(x_{\frac{j}{2^m}},\delta) \cap X$,

\begin{equation}
\labell{intInX}
\left[x_{\frac{j-1}{2^m}}, x_{\frac{j}{2^m}}\right] \subset X
\quad \text{for each} \quad  1 \leq j \leq 2^m.
\end{equation}
Since also $x_{\frac{j+1}{2^m}} \in B(x_{\frac{j}{2^m}},\delta)$,
\begin{equation*}
\left[x_{\frac{j-1}{2^m}}, x_{\frac{j+1}{2^m}}\right] \subset X
\quad \text{for each} \quad 1 \leq j < 2^m.
\end{equation*}
It follows that
\begin{equation*}
d_X(x_{\frac{j-1}{2^m}}, x_{\frac{j}{2^m}}) = \|x_{\frac{j-1}{2^m}} -
x_{\frac{j}{2^m}}\| \quad \text{ and } \quad d_X(x_{\frac{j-1}{2^m}},
x_{\frac{j+1}{2^m}}) = \|x_{\frac{j-1}{2^m}} - x_{\frac{j+1}{2^m}}\|.
\end{equation*}

Thus equation \eqref{iteratedmidpoint} can be rewritten as

\begin{equation*}
\|x_\frac{j-1}{2^m} - x_\frac{j}{2^m}\| = \|x_\frac{j}{2^m} -
x_\frac{j+1}{2^m}\| = \frac{1}{2}\|x_{\frac{j-1}{2^m}} -
x_{\frac{j+1}{2^m}}\|,
\end{equation*} 
which implies, by the triangle inequality, that the points 
$x_\frac{j-1}{2^m}$, $x_\frac{j}{2^m}$, $x_\frac{j+1}{2^m}$ 
are collinear.  
This and \eqref{intInX} imply that $[x_0, x_1] \subset X$.
\end{proof}

\section{Local and global convexity of maps}
\labell{sec:local-global}

The Tietze-Nakajima theorem involves subsets of $\R^n$.
We will now consider spaces with maps to $\R^n$
that are not necessarily inclusion maps.

Consider a continuous path $\gammaRn \colon [0,1] \to \R^n$.
Its length, which is denoted $l(\gammaRn)$, is the supremum,
over all natural numbers $N$ and all partitions 
$0 = t_0 < t_1 < \ldots < t_N = 1$, of 
$\sum_{i=1}^N \| \gammaRn(t_i) - \gammaRn(t_{i-1}) \|$.
We have $l(\gammaRn) \geq \| \gammaRn(1) - \gammaRn(0) \|$
with equality if and only if one of two cases occurs:
\begin{enumerate}
\item[(a)]
The path $\gammaRn$ is constant.
\item[(b)]
The image of $\gammaRn$ is the segment $[\gammaRn(0),\gammaRn(1)]$,
and $\gammaRn$ is a weakly monotone parametrization of this segment:
if $0 \leq t_1 < t_2 < t_3 \leq 1$, then the point $\gammaRn(t_2)$ 
lies on the segment $[\gammaRn(t_1),\gammaRn(t_3)]$.
\end{enumerate}

\begin{Definition}
The path $\gammaRn \colon [0,1] \to \R^n$ is
\emph{monotone straight} if it satisfies (a) or (b).
\end{Definition}

\begin{Definition} \labell{def:convex}
Let $X$ be a Hausdorff topological space.  A continuous map $\Psi$
from $X$ to $\R^n$, or to a subset of $\R^n$, is called \emph{convex} 
if every two points $x_0$ and $x_1$ in $X$ can be connected by a 
continuous path $\gamma \colon [0,1] \to X$ such that 
\begin{equation} \labell{condition}
\gamma(0) = x_0, \quad \gamma(1) = x_1, \quad 
\text{and} \quad \Psi \circ \gamma \text{ is monotone straight.}
\end{equation}
\end{Definition}

\begin{warning}
%
%
For a function $\psi$ from $\R$ to $\R$, the condition in 
Definition~\ref{def:convex} is equivalent to 
$\psi \colon \R \to \R$ being weakly monotone.
This is different from the usual notion of a convex function
(that $\psi(ta + (1-t)b) \leq t\psi(a) + (1-t)\psi(b)$ 
for all $a,b$ and for all $0 \leq t \leq 1$).
In the usual notion of a convex function, 
the domain $X$ must be an affine space,
and the target space must be $\R$.
In Definition~\ref{def:convex}, 
the domain $X$ is only a topological space,
and the target space can be $\R^n$.
In this paper, ``convex map" is always in the sense
of Definition~\ref{def:convex}.
\end{warning}

\begin{Remark} \labell{weakly monotone}
If $\Psi(x_0) = \Psi(x_1)$, condition \eqref{condition} 
means that the path $\gamma$ lies entirely within a level set of $\Psi$.
If $\Psi(x_0) \neq \Psi(x_1)$, the condition implies that the 
image of $\Psi \circ \gamma$ is the segment $[\Psi(x_0),\Psi(x_1)]$.
\end{Remark}


\begin{Example}
Consider the two-sphere $S^2 = \{ x \in \R^3 \ | \ \|x\|^2 = 1 \}$.
The height function $\Psi \colon S^2 \to \R$, given by 
$(x_1,x_2,x_3) \mapsto x_3$, is convex.  The projection
$\Psi \colon S^2 \to \R^2$, given by
$(x_1,x_2,x_3) \mapsto (x_1,x_2)$, is not convex.
\end{Example}

We shall prove the following generalization of the Tietze-Nakajima theorem:

\begin{Theorem} \labell{Theorem}
Let $X$ be a connected Hausdorff topological space, 
let $\calT \subset \R^n$ be a convex subset, and let
$$ \Psi \colon X \rightarrow \calT $$ 
be a continuous and proper map. 
Suppose that for every point $x \in X$ there exists an open 
neighbourhood $U \subset X$ of $x$ such that the map
$\Psi|_U \colon U \to \Psi(U)$
is convex and open.  Then
\begin{enumerate}
\item[(a)]
The image of $\Psi$ is convex.
\item[(b)]
The level sets of $\Psi$ are connected.
\item[(c)]
The map $\Psi \colon X \to \Psi(X)$ is open.
\end{enumerate}
\end{Theorem}

\begin{Remark} \labell{Tietze follows}
The Tietze-Nakajima theorem is the special case of Theorem \ref{Theorem}
in which $\calT = \R^n$, the space $X$ is a subset of~$\R^n$,
and the map $\Psi \colon X \to \R^n$ is the inclusion map.
\end{Remark}

The convexity of a map has the following immediate consequences:

\begin{Lemma} \labell{restrict}
If $\Psi \colon X \to \R^n$ is a convex map then,
for any convex subset $A \subset \R^n$,
the restriction of $\Psi$ to the preimage $\Psi\inv(A)$
is also a convex map.
\end{Lemma}

\begin{proof}
Let $A \subset \R^n$ be convex, and let $x_0, x_1 \in \Psi\inv(A)$.
Let $\gamma \colon [0,1] \rightarrow X$ be a path from $x_0$ to $x_1$
whose composition with $\Psi$ is monotone straight.
The image of $\Psi \circ \gamma$ is the (possibly degenerate) segment 
$[\Psi(x_0), \Psi(x_1)]$.
Because $A$ is convex and contains the endpoints of this segment,
it contains the entire segment, so $\gamma$ is a path in $\Psi\inv(A)$.
Thus, $x_0$ and $x_1$ are connected by a path in $\Psi\inv(A)$
whose composition with $\Psi$ is monotone straight.
\end{proof}


\begin{Lemma}[Global properties imply convexity]
\labell{properties imply convexity}
If $\Psi \colon X \to \R^n$ is a convex map, then its image, $\Psi(X)$, 
is convex, and its level sets, $\Psi\inv(w)$, for $w \in \Psi(X)$,
are connected.
\end{Lemma}

\begin{proof}
Take any two points in $\Psi(X)$; write them as 
$\Psi(x_0)$ and $\Psi(x_1)$ where $x_0$ and $x_1$ are in $X$.
Because the map $\Psi$ is convex,
there exists a path $\gamma$ in $X$ that connects
$x_0$ and $x_1$ and such that
the image of $\Psi \circ \gamma$ is the segment $[\Psi(x_0), \Psi(x_1)]$.
In particular, the segment $[\Psi(x_0),\Psi(x_1)]$ is contained 
in the image of $\Psi$.  This shows that the image of $\Psi$ is convex.

Now let $x_0$ and $x_1$ be any two points in $\Psi\inv(w)$.
Because the map $\Psi$ is convex, there exists a path $\gamma$ 
that connects $x_0$ and $x_1$ and such that
the curve $\Psi \circ \gamma$ is constant.
Thus, this curve is entirely contained in the level set $\Psi\inv(w)$.  
This shows that the level set $\Psi\inv(w)$ is connected.
\end{proof}

\begin{Remark} \labell{path lifting}
Suppose that the map $\Psi \colon X \to \R^n$ has the 
\emph{path lifting property}, i.e., for every path 
$\ol{\gamma} \colon [0,1] \to \R^n$
and every point $x \in \Psi\inv(\ol{\gamma}(0))$
there exists a path $\gamma \colon [0,1] \to X$
such that $\gamma(0) = x$ and $\Psi \circ \gamma = \ol{\gamma}$. 
Then the converse of Lemma~\ref{properties imply convexity}
holds: if the image $\Psi(X)$ is convex and the level sets
$\Psi\inv(w)$, $w \in \Psi(X)$, are path connected,
then the map $\Psi \colon X \to \R^n$ is convex.
\end{Remark}

The main ingredient in the proof of Theorem \ref{Theorem} 
is the following theorem, which we shall prove in section~\ref{sec:last proof}:

\begin{Theorem}
[Local convexity and openness imply global convexity and openness]
\labell{maintheorem}
Let $X$ be a connected Hausdorff topological space, 
let $\calT$ be a convex subset of $\R^n$, and let
$\Psi \colon X \to \calT$ be a continuous proper map.  
Suppose that for every point $x \in X$ 
there exists an open neighbourhood $U$ of $x$ such that the map
$ \Psi|_U \colon U \to \Psi(U)$
is convex and open.

Then the map 
$ \Psi \colon X \rightarrow \Psi(X)$
is convex and open.
\end{Theorem}

Following~\cite{HNP}, one may call Theorem \ref{maintheorem}
a \emph{Lokal-Global-Prinzip}.

\begin{Remark} \labell{U not small}
In Theorem \ref{maintheorem}, we assume that each point
is contained in an open set on which the map is convex
and is open as a map to its image, but we do not insist
that these open sets form a basis to the topology.
This requirement would be too restrictive, as is illustrated
in the following two examples.
\begin{enumerate}
\item[(a)]
Consider the map $(x,y) \mapsto -y + \sqrt{x^2 + y^2} $
from $\R^2$ to $\R$.  One level set is the non-negative
$y$-axis $\{ (0,y) \ | \ y \geq 0 \}$; the other level sets
are the parabolas $y = \frac{1}{2\alpha} x^2 - \frac{\alpha}{2}$
for $\alpha > 0$.  This map is convex, but its restrictions
to small neighborhoods of individual points on the positive
$y$-axis are not convex.  (These restrictions have
disconnected fibres.)
\item[(b)]
Consider the map $(t,e^{i\theta}) \mapsto t e^{i\theta}$
from $\R \times S^1$ to $\C \cong \R^2$.   This map is convex,
but its restrictions to small neighborhoods of individual
points on the zero section $\{ 0 \} \times S^1$ are not convex.
(These restrictions have a non-convex image.)
\end{enumerate}
\end{Remark}

\smallskip

\begin{proof} 
[Proof of Theorem \ref{Theorem}, assuming Theorem \ref{maintheorem}]
By Theorem \ref{maintheorem}, the map $\Psi$ is convex, and it is open
as a map to its image.
By Lemma \ref{properties imply convexity}, 
the level sets of $\Psi$ are connected and the image of $\Psi$ is convex.
\end{proof}

The bulk of this paper is devoted to proving Theorem \ref{maintheorem}.

\section{Convexity for components of preimages of neighbourhoods}

We first set some notation.

Let $X$ be a Hausdorff topological space
and $\Psi \colon X \to \R^n$ a continuous map.
For $x \in X$ with $\Psi(x) = w$, 
we denote by $[x]$ the path connected component of $x$ in $\Psi^{-1}(w)$, 
and, for $\eps > 0$, we denote by $U_{[x], \varepsilon}$ 
the path connected component of $x$
in $\Psi^{-1}(B(w, \varepsilon))$. Note that $U_{[x], \varepsilon}$
does not depend on the particular choice of $x$ in $[x]$.

\begin{remark}
Suppose that every point in $X$ has an open neighbourhood $U$
on which the restriction $\Psi|_U$ is convex.
Then, in the definitions of $[x]$ and $U_{[x],\vareps}$,
the term \emph{path connected component} can be replaced by 
\emph{connected component}.
Indeed, let $Y = \Psi\inv(B(w,\eps))$ or $Y = \Psi\inv(w)$.
If $\Psi|_U$ is convex, so is $\Psi|_{U \cap Y}$;
in particular, $U \cap Y$ is path connected.
Thus, every point in $Y$ has a path connected
open neighborhood with respect to the relative topology on $Y$. 
So the connected components of $Y$ coincide with its
path connected components.
\end{remark}

A crucial step in the proof of Theorem \ref{maintheorem}
is that the neighbourhoods $U$ such that 
$\Psi|_U \colon U \to \Psi(U)$
is convex and open can be taken to be the 
entire connected components $U_{[x],\vareps}$:

\begin{prop}[Properties for connected components]
\labell{intgoalB}
Let $X$ be a Hausdorff topological space, 
$\calT \subset \R^n$ a convex subset, and
$\Psi \colon X \to \calT$ a continuous proper map.  
Suppose that for every point $x \in X$ 
there exists an open neighbourhood $U$ of $x$ such that the map
$\Psi|_U \colon U \to \Psi(U)$ is convex and open.

Then for every point $x \in X$ there exists an $\vareps > 0$
such that the map
$\Psi|_{U_{[x],\vareps}} \colon U_{[x],\vareps} \to \Psi(U_{[x],\vareps})$
is convex and open.
\end{prop}

We digress to recall standard consequences of the properness of a map.

\begin{Lemma} \labell{consequences of proper}
Let $X$ be a Hausdorff topological space, $\calT \subset \R^n$
a subset, and $\Psi \colon X \to \calT$ a continuous proper map.
Let $w_0 \in \calT$.
\begin{enumerate}
\item
Let $U$ be an open subset of $X$ that contains 
the level set $\Psi\inv(w_0)$.  Then there exists $\eps > 0$ 
such that the pre-image $\Psi\inv(B(w_0,\eps))$ is contained in $U$.

\item
Suppose that every point of $\Psi\inv(w_0)$
has a connected open neighborhood in $\Psi\inv(w_0)$
with respect to the relative topology.
Then there exists $\eps > 0$ such that 
whenever $[x]$ and $[y]$ are distinct connected components
of $\Psi\inv(w_0)$
the sets $U_{[x],\eps}$ and $U_{[y],\eps}$ are disjoint.
\end{enumerate}
\end{Lemma}

\begin{proof}[Proof of part (1)]
Suppose otherwise.  Then, for every $\varepsilon > 0$ there
exists $x_\varepsilon \in X \smallsetminus U$ 
such that $\| \Psi(x_\varepsilon) - w_0 \| < \varepsilon$.

Let $\varepsilon_j$ be a sequence such that 
$\varepsilon_j \rightarrow 0$ as $j \rightarrow \infty$.  
Then $x_{\varepsilon_j} \in X \smallsetminus U$ for all $j$, 
and $\Psi(x_{\varepsilon_j}) \rightarrow w_0$ as $j \rightarrow \infty$.

The set 
$\{\Psi(x_{\varepsilon_j})\}_{j=1}^\infty \cup \{w_0\}$ 
is compact.  By properness, its preimage,
$\cup_{j=1}^\infty \Psi^{-1}(\Psi(x_{\varepsilon_j})) \cup \Psi^{-1}(w_0)$,
is compact.  
The sequence $\{x_{\varepsilon_j}\}_{j=1}^\infty$ is in this preimage.  
So there exists a point $x_\infty$ such that
every neighborhood of $x_\infty$ contains $x_{\eps_j}$ 
for infinitely many values of $j$.

By continuity, $\Psi(x_\infty) = w_0$.  
Since $U$ contains $\Psi^{-1}(w_0)$ and is open,
%
%
$U$ is a neighborhood of $x_\infty$, so there exist arbitrarily 
large values of $j$ such that $x_{\varepsilon_j} \in U$. 
This contradicts the assumption $x_{\eps_{j}} \in X \ssminus U$.
\end{proof}

\begin{proof}[Proof of part (2)]
Because $\Psi$ is proper, the level set $\Psi\inv(w_0)$ is compact.  
Because $\Psi\inv(w_0)$ is compact and is covered by connected
open subsets with respect to the relative topology,
it has only finitely many components $[x]$.
Because these components are compact and disjoint and $X$ is Hausdorff,
there exist open subsets $\mathcal{O}_{[x]}$ in $X$
such that $[x] \subset \mathcal{O}_{[x]}$ for each component $[x]$
of $\Psi\inv(w_0)$
and such that for $[x]$ and $[y]$ in $\Psi^{-1}(w_0)$,
if $[x] \neq [y]$ then $\mathcal{O}_{[x]} \cap \mathcal{O}_{[y]} = \emptyset$.
The union of the sets $\mathcal{O}_{[x]}$ is an open subset of $X$
that contains the fiber $\Psi\inv(w_0)$.
By part (1), this open subset contains $\Psi\inv(B(w_0,\vareps))$ 
for every sufficiently small $\eps$.
For such an~$\vareps$, 
because each $U_{[x],\vareps}$ is contained in $\mathcal{O}_{[x]}$
and the sets $\mathcal{O}_{[x]}$ are disjoint,
the sets $U_{[x],\vareps}$ are disjoint.
\end{proof}

We now prepare for the proof of Proposition~\ref{intgoalB}.
In the remainder of this section,
let $X$ be a Hausdorff topological space, $\calT \subset \R^n$
a subset, and $\Psi \colon X \to \calT$ a continuous map.
Fix a point $w_0 \in \calT$.
Let $\{ U_i \}$ be a collection of open subsets of $X$
whose union contains $\Psi\inv(w_0)$.



\begin{Lemma} \labell{sequence}
Let $[x]$ be a connected component of $\Psi\inv(w_0)$.
If $U_k \cap [x] \neq \emptyset$ and $U_l \cap [x] \neq \emptyset$,
then there exists a sequence $k= i_0, i_1, \ldots, i_s=l$
%
such that
\begin{equation} \labell{adjacent}
 U_{i_{q-1}} \cap U_{i_q} \cap [x] \neq \emptyset
\quad \text{ for } \quad q=1,\ldots,s.
\end{equation}
\end{Lemma}

\begin{proof}
Let $I_k$ denote the set of indices $j$ for which one can get
from $U_k$ to $U_j$ through a sequence of sets 
%
%
$U_{i_0}, U_{i_1}, \ldots, U_{i_s}$ with the property~\eqref{adjacent}.
If $j \in I_k$ and $j' \not \in I_k$
then $U_j \cap [x]$ and $U_{j'} \cap [x]$ are disjoint.  Thus
\begin{displaymath}
[x] = \left( \bigcup_{j \in I_k}      U_j \cap [x] \right) 
 \cup \left( \bigcup_{j' \not\in I_k}  U_{j'} \cap [x] \right)
\end{displaymath}
expresses $[x]$ as the union of two disjoint open subsets,
of which the first is non-empty.  Because $[x]$ is connected,
the second set in this union must be empty.
So $U_l \cap [x] \neq \emptyset$ implies $l \in I_k$.
\end{proof}


Now assume, additionally, that the covering $\{ U_i \}$ is finite
and that, for each $i$, the map 
$\Psi|_{U_i} \colon U_i \to \Psi(U_i)$
is open.  Let 
\begin{equation} \labell{def of Wi}
W_i := \Psi(U_i).
\end{equation}

\begin{Lemma} \labell{closelythesame}
Let $[x]$ be a connected component of $\Psi\inv(w_0)$.
For sufficiently small $\varepsilon > 0$, the following is true.
\begin{enumerate}
\item 
For any $i$ and $j$, 
%
%
if $U_i \cap U_j \cap [x] \neq \emptyset$, then
\begin{displaymath}
W_i \cap B(w_0, \varepsilon) = \Psi(U_i \cap U_j) \cap B(w_0, \varepsilon) .
\end{displaymath}

\item
For any $k$ and $l$, 
%
%
if $U_k \cap [x]$ and $U_l \cap [x]$ are non-empty, then
\begin{displaymath}
W_k \cap B(w_0, \varepsilon) = W_l \cap B(w_0, \varepsilon) .
\end{displaymath}
\end{enumerate}
\end{Lemma}


\begin{proof}
%
%
Suppose that $U_i \cap U_j \cap [x] \neq \emptyset$.
Then the set $\Psi(U_i \cap U_j)$ contains $w_0$.  
Since $U_i \cap U_j$ is open in $U_i$, 
and since the restriction of $\Psi$ to $U_i$ is an open map to its image,
the set $\Psi(U_i \cap U_j)$ is open in $W_i$.  
Let $\varepsilon_{ij} > 0$ be such that the set $\Psi(U_i \cap U_j)$
contains $W_i \cap B(w_0, \varepsilon_{ij})$.
Because we also have $\Psi(U_i \cap U_j) \subset \Psi(U_i) = W_i$,
\begin{displaymath}
W_i \cap B(w_0, \varepsilon_{ij}) 
 = \Psi(U_i \cap U_j) \cap B(w_0, \varepsilon_{ij}).
\end{displaymath}

Let $\vareps$ be any positive number that is smaller than
or equal to $\vareps_{ij}$ 
%
%
for all the pairs $U_i$, $U_j$ for which $U_i \cap U_j \cap [x] 
\neq \emptyset$.  Then, for every such pair $U_i$, $U_j$,
$$ W_i \cap B(w_0,\vareps) = \Psi(U_i \cap U_j) \cap B(w_0,\vareps).$$
This proves (1).

Now suppose that $U_k \cap [x] \neq \emptyset$ 
and $U_l \cap [x] \neq \emptyset$.  By Lemma \ref{sequence},
one can get from $U_k$ to $U_l$ by a sequence of sets 
%
%
$U_k = U_{i_0}, \ldots U_{i_s} = U_l$ such that
$U_{i_{q-1}} \cap U_{i_q} \cap [x] \neq \emptyset$ for $q=1,\ldots,s$.
Part (1) then implies that the intersections $W_{i_q} \cap B(w_0,\eps)$
are the same for all the elements in the sequence.
Because the sequence
begins with $U_k$ and ends with $U_l$, it follows that
$$ W_k \cap B(w_0,\vareps) = W_l \cap B(w_0,\vareps).$$
This proves (2).
\end{proof}


Let $[x]$ be a connected component of $\Psi\inv(w_0)$.
Fix an $\varepsilon > 0$ that satisfies the conditions of
Lemma \ref{closelythesame}.  
%
%
Let
\begin{equation}
\labell{Wx}
\Wxeps := W_i \cap B(w_0, \varepsilon) 
\qquad \text{ when } U_i \cap [x] \neq \emptyset.
\end{equation} 
By part (2) of Lemma \ref{closelythesame}, 
this is independent of the choice of such $i$.  
Also, define 
\begin{equation} \labell{Uxeps}
\tU_{[x],\vareps} := \bigcup_{\substack{U_i \cap [x] \neq \emptyset}}
         U_i \cap \Psi^{-1}(B(w_0, \varepsilon)).
\end{equation}

We have
\begin{equation} \labell{image is Wi}
\begin{aligned}
\Psi(\tUxeps) &= \bigcup_{\substack{U_i \cap [x] \neq \emptyset}}
   \Psi(U_i) \cap B(w_0, \varepsilon) 
                                 \qquad \text{by \eqref{Uxeps}} \\
 &=  W_{[x], \varepsilon} \qquad \text{by \eqref{def of Wi} and \eqref{Wx}} .
\end{aligned}
\end{equation}


\begin{Lemma} \labell{components1}
Suppose that, for each $i$, the level sets of 
$ \Psi|_{U_i} \colon U_i \to W_i $
are path connected.
Then
%
%
the level sets of
$ \Psi|_{\tU_{[x],\vareps}} \colon 
         \tU_{[x],\vareps} \to W_{[x],\vareps} $
are path connected.
\end{Lemma}

\begin{proof}
%
%
%
%

Let $w \in \Wxeps$ and let $x_0, x_1 \in \tUxeps \cap \Psi^{-1}(w)$. 
By \eqref{Uxeps} there exist $i$ and $k$ such that 
$x_0 \in U_i$, $x_1 \in U_k$, $U_i \cap [x] \neq \emptyset$, 
and $U_k \cap [x] \neq \emptyset$.
Fix such $i$ and $k$.
%
%
By Lemma \ref{sequence}, 
there exists a sequence $i = i_0, i_1, \ldots, i_s = k$ 
such that $U_{i_{l-1}} \cap U_{i_l} \cap [x] \neq \emptyset$
for $l = 1, \ldots, s$.
Part (1) of Lemma~\ref{closelythesame},
and the definition~\eqref{Wx} of $\Wxeps$,
imply that $\Psi(U_{i_{l-1}} \cap U_{i_l} ) \cap B(w_0,\eps) = \Wxeps$,
and thus
$U_{i_{l-1}} \cap U_{i_l} \cap \Psi^{-1}(w)$
is non-empty for each $1 \leq l \leq s$.
Since each $U_{i_l} \cap \Psi^{-1}(w)$ is path connected,
this implies that $x_0$ and $x_1$ can be connected by a path 
in $\tUxeps \cap \Psi^{-1}(w)$.
\end{proof}

\begin{Lemma} \labell{components2}
Suppose that, for each $i$, the restriction of $\Psi$ to $U_i$
is a convex map. Then the map 
\begin{equation} \labell{rest}
 \Psi|_{\tUxeps} \colon \tUxeps \to \Psi(\tUxeps) 
\end{equation}
is convex and open.
\end{Lemma}

\begin{proof}
Let $x_0$ and $x_1$ be in $\tUxeps$.
Let $i$ be such that $x_0 \in U_i$ and $U_i \cap [x] \neq \emptyset$.
By \eqref{image is Wi}, $\Psi(x_1) \in \Wxeps$.
By \eqref{Wx} and \eqref{def of Wi},
there exists $y \in U_i$ such that $\Psi(y) = \Psi(x_1)$.

By assumption, the restriction of~$\Psi$ to~$U_i$ is a convex map.
By Lemma~\ref{restrict}, the restriction of~$\Psi$
to~$U_i \cap \Psi\inv(B(w_0,\eps))$ is also convex.
Let $\gamma'$ be a path in $U_i \cap \Psi\inv(B(w_0,\eps))$ 
from $x_0$ to $y$ such that $\psi \circ \gamma'$ is monotone straight.
By Lemma \ref{components1}
there exists a path $\gamma''$ in $\tU_{[x],\vareps}$
from $y$ to $x_1$
whose composition with $\Psi$ is constant.
Let $\gamma$ be the concatenation of $\gamma'$ with $\gamma''$;
then $\gamma$ is a path from $x_0$ to $x_1$ 
and $\Psi \circ \gamma$ is monotone straight.

Thus, the map \eqref{rest} is convex.
To show that this map is open, we want to show
that given any open set $\Omega \subset \tUxeps$, 
its image $\Psi(\Omega)$ is open in $\Wxeps$.  By \eqref{Uxeps}, 
$\Psi(\Omega) = \cup_i \Psi(\Omega \cap U_i)$ 
for $i$ such that $U_i \cap [x] \neq \emptyset$, 
and each $\Psi(\Omega \cap U_i)$ is contained in $B(w_0, \varepsilon)$.  
Since $\Psi|_{U_i} \colon U_i \rightarrow W_i$ is open, 
$\Psi(\Omega \cap U_i)$ is open in $W_i$.  
By \eqref{Wx}, each $\Psi(\Omega \cap U_i)$ is open in $\Wxeps$.
\end{proof}

\begin{proof}[Proof of Proposition \ref{intgoalB}]
Let $w_0 = \Psi(x)$.  For each $x' \in \Psi\inv(w_0)$,
let $U_{x'}$ be an open neighbourhood of $x'$
such that the map $\Psi|_{U_{x'}} \colon U_{x'} \to \Psi(U_{x'})$
is convex and open.
The sets $U_{x'}$, for $x' \in \Psi\inv(w_0)$, 
cover $\Psi\inv(w_0)$.
Because $\Psi\inv(w_0)$ is compact, there exists
a finite subcovering; let $\{U_i\}_{i=1}^n$ be a finite subcovering.

Because $\Psi\inv(w_0)$ is compact and each point has a connected
neighborhood with respect to the relative topology,
$\Psi\inv(w_0)$ has only finitely many components $[x]$.
%
%
Let $\eps > 0$ satisfy the conditions of Lemma~\ref{closelythesame}
for all these components.
By Lemma~\ref{consequences of proper},
after possibly shrinking $\eps$, we may assume that
$\Psi\inv(B(w_0,\eps)) \subset \cup_i U_i$
and that whenever $[x]$ and $[y]$
are distinct connected components of $\Psi\inv(w_0)$
the sets $\Uxeps$ and $\Uyeps$ are disjoint.

Let $\tUxeps$ and $\Wxeps$ be the sets defined 
in~\eqref{Uxeps} and~\eqref{Wx}.
Then the preimage $\Psi\inv(B(w_0,\eps))$ is the union of the sets
$\tUxeps$, for components $[x]$ of $\Psi\inv(w_0)$.
Because each $\tUxeps$ is connected and contains $[x]$, it is contained
in the connected component $\Uxeps$ of $x$ in $\Psi\inv(B(w_0,\eps))$.
Because the sets $\Uxeps$ are disjoint and the union 
of the sets $\tUxeps$ is the entire preimage $\Psi\inv(B(w_0,\eps))$,
each $\tUxeps$ is \emph{equal} to $\Uxeps$.
This and Lemma~\ref{components2} give Proposition~\ref{intgoalB}.
\end{proof}

\section{Distance with respect to a locally convex map}

Let $X$ be a Hausdorff topological space and $\Psi \colon X \to \R^n$
a continuous map.  Let $x_0$ and $x_1$ be two points in $X$. 
We define their $\Psi$-distance to be
\begin{equation*} 
d_{\Psi}(x_0,x_1) = 
   \inf\{l(\Psi \circ \gamma) \ | \ 
         \gamma \colon \left[0,1\right] \rightarrow X, \ 
         \gamma(0)=x_0, \ \gamma(1)=x_1 \}.
\end{equation*}
Note that the $\Psi$-distance can take any value in $[0,\infty]$.
Also note that $d_\Psi(x_0,x_1) = 0$ if and only if $x_0$ and $x_1$ 
are in the same path-component of a level set of $\Psi$.


\begin{remark}
In practice, we will work with a space $X$ which is connected
and in which each point has a neighbourhood $U$ such that
the restriction of $\Psi$ to $U$ is a convex map.
For such a space, in the above definition of $\Psi$-distance,
we may take the infimum to be over the set of paths $\gamma$
such that $\Psi \circ \gamma$ is polygonal:

Indeed, let $\gamma \colon [0,1] \rightarrow X$ be any path 
such that $\gamma(0) = x_0$ and $\gamma(1) = x_1$.

By our assumption on $X$,
for every $\tau \in [0,1]$ there exists an open interval $J_\tau$
containing $\tau$ and an open subset $U_\tau \subset X$
such that the restriction of $\Psi$ to $U_\tau$ is a convex map
and such that $\gamma(J_\tau \cap [0,1]) \subset U_\tau$.

The open intervals $\{ J_\tau \}$ form an open covering of $[0,1]$.
Because the interval $[0,1]$ is compact, there exists a finite
subcovering; denote it $J_1,\ldots,J_s$.  Let
$$\eps = \min \{ \, \text{length}(J_i \cap J_k) | \, 
i,k \in \{1,\ldots,s\} \text{ and } J_i \cap J_k \neq \emptyset \}.$$

Any subinterval $[\alpha,\beta] \subset [0,1]$
of length $< \eps$ is contained in one of the $J_i$s.
Indeed, given such a subinterval $[\alpha,\beta]$,
consider those intervals of $J_1,\ldots,J_s$ that contain $\alpha$; 
let $J_i$ be the one whose upper bound $b_i$,
is maximal; then $J_i$ also contains $\beta$.

Thus, for any subinterval $[\alpha,\beta] \subset [0,1]$
of length $< \eps$ there exists an open subset $U \subset X$
such that the restriction of $\Psi$ to $U$ is a convex map
and such that $\gamma(\alpha)$ and $\gamma(\beta)$ are both 
contained in $U$.

Partition $[0,1]$ into $m$ intervals $0 = t_0 < \ldots < t_m = 1$ such that 
$|t_j - t_{j-1}| < \varepsilon$ for each $j$.
By the previous paragraph, for every $1 \leq j \leq m$
there exists $U \subset X$ 
such that the restriction of $\Psi$ to $U$ is a convex map
and such that $\gamma(t_{j-1})$ and $\gamma(t_j)$ are both
contained in $U$. Because the restriction of $\Psi$ to $U$ is convex,
there exists a path $\gamma_j$ 
in $X$ connecting $\gamma(t_{j-1})$ and $\gamma(t_j)$ such that 
the image of $\Psi \circ \gamma$ is a (possibly degenerate) segment
with a weakly monotone parametrization.
The path $\gamma'$ that is formed by concatenating
$\gamma_1, \ldots, \gamma_m$ connects $x_0$ and $x_1$, 
the composition $\Psi \circ \gamma'$ is polygonal,
and $l(\Psi\circ\gamma') \leq l(\Psi\circ\gamma)$.
\end{remark}


\section{Proof that local convexity and openness
         imply global convexity and openness}
\labell{sec:last proof}

\begin{Lemma}[Existence of midpoint] \labell{Midpoint}
Let $X$ be a connected Hausdorff topological space,
$\calT \subset \R^n$ a convex subset, and
$\Psi \colon X \rightarrow \calT$ a continuous and proper map.
Suppose that for every point $x \in X$ there exists 
an open neighbourhood $U$ such that the restriction of $\Psi$
to $U$ is a convex map.

Let $x_0$ and $x_1$ be in $X$.  Then there exists a point $x_{1/2} \in X$ 
such that
\begin{equation}
\labell{xhalf}
d_{\Psi}(x_0, x_{1/2}) = 
   d_{\Psi}(x_{1/2}, x_1) = \frac{1}{2}d_{\Psi}(x_0, x_1).
\end{equation}
\end{Lemma}

\begin{proof}
Choose paths $\gamma_n$ connecting $x_0$ and $x_1$ such that
the sequence
$\{l(\Psi \circ \gamma_n)\}$ converges to $d_{\Psi}(x_0, x_1)$.

Let $t_j \in [0,1]$ be such that $\gamma_j(t_j)$ is the midpoint of
the path $\gamma_j$:

\begin{equation*}
l(\Psi \circ \gamma_j \arrowvert_{\scriptscriptstyle{\left[0,t_j\right]}}) = 
l(\Psi \circ \gamma_j \arrowvert_{\scriptscriptstyle{\left[t_j,1\right]}}) =
\frac{1}{2}l(\Psi \circ \gamma_j).
\end{equation*} 

Let $r > \half d_{\Psi}(x_0,x_1)$. 
Then all but finitely many of the midpoints $\gamma_j(t_j)$
lie in the set
\begin{equation*}
A = \Psi^{-1} ( \overline{B}(\Psi(x_0),r) ) .
\end{equation*}
This set is compact because $\Psi$ is proper.
So there exists a point $\xhalf$ such that every neighbourhood
of $\xhalf$ contains $\gamma_j(t_j)$ for infinitely many values of $j$.
We will show that the point $\xhalf$ satisfies equation~\eqref{xhalf}.

We first show that
$d_{\Psi}(x_0,x_\frac{1}{2}) \leq \frac{1}{2}d_{\Psi}(x_0,x_1)$, or,
equivalently, that for every $\varepsilon > 0$ there exists a path
$\gamma$ connecting $x_0$ and $x_\frac{1}{2}$ such that 
$l(\Psi \circ \gamma) < \frac{1}{2}d_{\Psi}(x_0,x_1) + \varepsilon$.

Let $U$ be a neighbourhood of $x_{1/2}$ such that the restriction
of $\Psi$ to $U$ is a convex map.
Let $j$ be such that the following facts are true:

\begin{enumerate}
\item[(i)] $\gamma_j(t_j) \in U$ 
       and $\| \Psi(\gamma_{j}(t_{j})) - \Psi(\xhalf) \| 
            < \frac{\varepsilon}{2}$.
\item[(ii)] $l(\Psi \circ \gamma_j) < d_{\Psi}(x_0,x_1) + \varepsilon$.
\end{enumerate}

By (i) and since $\Psi|_U$ is a convex map,
there exists a path $\mu$
connecting $\gamma_j(t_j)$ and $\xhalf$ such that 
$l(\Psi \circ \mu) < \frac{\varepsilon}{2}$.  Let $\gamma$ be the
concatenation of $\gamma_j|_{[0,t_j]}$ and $\mu$.  Then $\gamma$ is a
path connecting $x_0$ and $\xhalf$, and 
$l(\Psi \circ \gamma) < \frac{1}{2}d_{\Psi}(x_0, x_1) + \varepsilon$.

Thus, $d_{\Psi}(x_0,x_{1/2}) \leq \frac{1}{2}d_{\Psi}(x_0,x_1)$.  By the
same argument, $d_{\Psi}(x_{1/2},x_1) \leq \frac{1}{2}d_{\Psi}(x_0,x_1)$.
If either of these were a strict inequality, then it would be
possible to construct a path from $x_0$ to $x_1$ whose image has length
less than $d_{\Psi}(x_0, x_1)$, 
which contradicts the definition of $d_\psi(x_0,x_1)$.  Thus,
$d_{\Psi}(x_0, x_{1/2}) = d_{\Psi}(x_{1/2}, x_1) =
\frac{1}{2}d_{\Psi}(x_0, x_1)$.  
\end{proof}

To prove Theorem \ref{maintheorem}, we need to have some uniform control 
on the sizes of $\vareps$ such that the restrictions of $\Psi$
to the connected components $U_{[x],\vareps}$ of $\Psi\inv(B(w_0,\eps))$ 
are convex.  The precise result that we will use is established 
in the following proposition:

\begin{Proposition} \labell{intgoalA}
Let $X$ be a Hausdorff topological space
and let $\Psi \colon X \to \R^n$ be a continuous map.
Suppose that for each $x \in X$ there exists an $\varepsilon > 0$ 
such that the restriction of $\Psi$ to the set $U_{[x],\vareps}$ 
is a convex map.

Then for every compact subset $A \subset X$ there exists 
$\varepsilon > 0$ such that for every $x \in A$ and $x' \in X$,
if $d_{\Psi}(x,x') < \varepsilon$, then there exists
a path $\gamma \colon [0,1] \to X$ such that $\gamma(0) = x$,
$\gamma(1) = x'$, and $\Psi \circ \gamma$ is monotone straight.
\end{Proposition}

\begin{proof}
For each $x \in X$, let $\vareps_x > 0$ be
such that the restriction of $\Psi$ to the set $U_{[x],\vareps_x}$ 
is a convex map.  The sets $U_{[x],\vareps_x/2}$,
for $x \in A$, form an open covering of the compact set $A$.
Choose a finite subcovering:  let $x_1,\ldots,x_k$ be points of $A$
and $\vareps_1,\ldots,\vareps_k$ be positive numbers such that,
for each $1 \leq i \leq k$, the restriction of $\Psi$
to the set $U_{[x_i],\vareps_i}$ is a convex map,
and such that the sets $U_{[x_i],\vareps_i/2}$ cover $A$.

Let $\vareps = \min\limits_{1 \leq i \leq k } \frac{\vareps_i}{2} $.

Let $x \in A$.
Let $1 \leq i \leq k$ be such that $x \in U_{[x_i],\vareps_i/2}$.

Because $U_{[x_i],\vareps_i/2}$, by its definition, is contained in
$\Psi\inv(B(\Psi(x_i),\vareps_i/2))$, we have 
$\| \Psi(x) - \Psi(x_i) \| < \vareps_i/2$.

Because $x$ and $x_i$ are also contained in the larger set 
$U_{[x_i],\vareps_i}$, and the restriction of $\Psi$ to this set
is a convex map, there exists a path $\gamma'$ from $x_i$ to $x$ 
such that $\Psi \circ \gamma'$ is monotone straight;
in particular,
$l(\Psi \circ \gamma') = \| \Psi(x) - \Psi(x_i) \|$,
so $l(\Psi \circ \gamma') < \vareps_i / 2$.

Let $x' \in X$ be such that $d_\Psi(x,x') < \vareps$.
Then, by the definition of $d_\Psi$,
there exists a path $\gamma''$ from $x$ to $x'$
such that $l(\Psi \circ \gamma'') < \vareps$. 

Let $\hat{\gamma}$ be the concatenation of $\gamma'$ and $\gamma''$.
Then $\hat{\gamma}$ is a path from $x_i$ to $x'$, and
$l(\Psi \circ \hat{\gamma}) \leq l(\Psi \circ \gamma') + l(\Psi \circ \gamma'')
 < \vareps_i/2 + \vareps \leq \vareps_i$.
Therefore, $\Psi \circ \hat{\gamma} \subset B(\Psi(x_i),\vareps_i)$.
Thus, $x'$ and $x_i$ are in the same connected component
of $\Psi\inv(B(\Psi(x_i),\vareps_i))$; that is,
$x'$ is in the set $U_{[x_i],\vareps_i}$. 
Because $x$ is also in the set $U_{[x_i],\vareps_i}$,
and because the restriction of $\Psi$ to this set is a convex map,
there exists a path $\gamma$ from $x$ to $x'$ 
such that $\Psi \circ \gamma$ is monotone straight.
\end{proof}


\begin{Proposition} \labell{tietzepsi}
Let $X$ be a connected Hausdorff topological space,
let $\calT \subset \R^n$ be a convex subset, 
and let $\Psi \colon X \to \calT$ be a continuous proper map.
Suppose that for every compact subset $A \subset X$ there exists 
$\varepsilon > 0$ such that, for every $x \in A$ and $x' \in X$,
if $d_{\Psi}(x,x') < \varepsilon$, then there exists
a path $\gamma \colon [0,1] \to X$ such that $\gamma(0) = x$,
$\gamma(1) = x'$, and $\Psi \circ \gamma$ is monotone straight.

Then $\Psi \colon X \to \R^n$ is a convex map.
\end{Proposition}

\begin{proof}
Fix $x_0$ and $x_1$ in $X$.  

By Lemma \ref{Midpoint}, there exists a point $x_{1/2}$ such that
\begin{displaymath}
d_{\Psi}(x_0, x_{1/2}) 
   = d_{\Psi}(x_{1/2}, x_1) 
   = \frac{1}{2}d_{\Psi}(x_0, x_1).
\end{displaymath}

Likewise, there exists a point $x_{1/4}$ which satisfies
\begin{displaymath}
d_{\Psi}(x_0, x_{1/4}) 
  = d_{\Psi}(x_{1/4}, x_{1/2}) 
  = \frac{1}{2}d_{\Psi}(x_0, x_{1/2}).
\end{displaymath}

By iteration, we get a map $\frac{j}{2^m} \mapsto x_{\frac{j}{2^m}}$,
for nonnegative integers $j$ and $m$ with $0 \leq j \leq 2^m$, 
such that
\begin{equation} \labell{eq1} 
d_{\Psi}(x_\frac{j-1}{2^m}, x_\frac{j}{2^m}) = d_{\Psi}(x_\frac{j}{2^m},
x_\frac{j+1}{2^m}) = \frac{1}{2}d_{\Psi}(x_\frac{j-1}{2^m},
x_\frac{j+1}{2^m}).
\end{equation}

Let $r > d_{\Psi}(x_0, x_1)$.
Let $\vareps > 0$ be associated with the compact set
\begin{displaymath}
A = \Psi^{-1}(\overline{B}(\Psi(x_0), r))
\end{displaymath}
as in the assumption of the proposition.

Choose $m$ large enough such that for every $1 \leq j \leq 2^m$,
\begin{displaymath}
d_{\Psi}(x_{\frac{j-1}{2^m}}, x_{\frac{j}{2^m}}) < \frac{\varepsilon}{2}.
\end{displaymath}

By the assumption, there exists a path $\gamma_j$ 
from $x_{(j-1)/{2^m}}$ to $x_{{j}/{2^m}}$ such that 
$\Psi \circ \gamma_j$ is monotone straight.
Thus,
\begin{equation} \labell{colin}
d_\Psi(x_{\frac{j-1}{2^m}}, x_{\frac{j}{2^m}}) =
\|\Psi(x_{\frac{j-1}{2^m}}) - \Psi(x_{\frac{j}{2^m}})\|
\quad \text{for each }  1 \leq j \leq 2^m .
\end{equation}

Similarly, 
$$ d_\Psi(x_{\frac{j-1}{2^m}}, x_{\frac{j+1}{2^m}}) =
   \|\Psi(x_{\frac{j-1}{2^m}}) - \Psi(x_{\frac{j+1}{2^m}})\| 
\quad \text{for each } 1 \leq j < 2^m.  $$

Thus equation \eqref{eq1} can be rewritten as
\begin{equation*}
\|\Psi(x_\frac{j-1}{2^m}) - \Psi(x_\frac{j}{2^m})\| 
= \|\Psi(x_\frac{j}{2^m}) - \Psi(x_\frac{j+1}{2^m})\| 
= \frac{1}{2}\|\Psi(x_{\frac{j-1}{2^m}}) - \Psi(x_{\frac{j+1}{2^m}})\|,
\end{equation*} 
which implies, by the triangle inequality, 
that the points $\Psi(x_\frac{j-1}{2^m})$,
$\Psi(x_\frac{j}{2^m})$, $\Psi(x_\frac{j+1}{2^m})$
are collinear.
The concatenation of the paths $\gamma_j$,
for $1 \leq j \leq 2^m$, is a path from $x_0$ to $x_1$
whose composition with $\Psi$ is monotone straight.
\end{proof}

\bigskip

\begin{proof}[Proof of Theorem \ref{maintheorem}]
Let $X$ be a connected Hausdorff topological space,
$\calT \subset \R^n$ a convex subset, and $\Psi \colon X \to \calT$
a continuous proper map.
Suppose that for every point $x \in X$
there exists an open neighbourhood $U$ such that the map
$\Psi|_U \colon U \to \Psi(U)$ is convex and open.

By Proposition \ref{intgoalB}, for every point $x$
there exists an $\vareps > 0$ such that the map
$$ \Psi|_{U_{[x],\vareps}} \colon U_{[x],\vareps} \to \Psi(U_{[x],\vareps})$$
is convex and open.

By Proposition \ref{intgoalA}, for every compact subset $A \subset X$
there exists an $\vareps > 0$ such that for every $x \in A$ and $x' \in X$,
if $d_\Psi(x,x') < \vareps$, then there exists a path 
$\gamma \colon [0,1] \to X$ such that $\gamma(0) = x$, $\gamma(1) = x'$,
and $\Psi \circ \gamma$ is monotone straight.

By Proposition \ref{tietzepsi}, the map $\Psi \colon X \to \R^n$
is convex.

\medskip

To show that the map $\Psi \colon X \to \Psi(X)$ is open,
it is enough to show that for each $w_0 \in \R^n$
there exists $\vareps > 0$ such that
the restriction of $\Psi$ to $\Psi\inv(B(w_0,\vareps))$
is open as a map to its image.

Fix $w_0 \in \R^n$.

Because the map $\Psi \colon X \to \R^n$ is convex,
the level set $\Psi\inv(w_0)$ is connected.
Thus, this level set consists of a single connected component, $[x]$.

By Proposition~\ref{intgoalB}, for sufficiently small $\vareps$,
the restriction of $\Psi$ to the set $U_{[x],\vareps}$ is open
as a map to its image.
The set $U_{[x],\vareps}$ is an open set that contains $\Psi\inv(w_0)$.
Because $\Psi$ is proper, there exists an $\vareps' > 0$
such that the set $U_{[x],\vareps}$ contains the preimage
$\Psi\inv(B(w_0,\vareps'))$; see Lemma \ref{consequences of proper}.
Thus, the restriction of $\Psi$ to the preimage $\Psi\inv(B(w_0,\vareps'))$
is open as a map to its image.
\end{proof}

\section{Examples}
\labell{sec:moment}

\begin{noTitle} \labell{Cn}
The map from $\C^n$ to $\R^n$ given by
\begin{equation} \labell{PhiCn}
 (z_1,\ldots,z_n) \mapsto (|z_1|^2,\ldots,|z_n|^2) 
\end{equation}
is convex, and it is open as a map from $\C^n$ to the
positive orthant $\R_+^n$.

Moreover, the restriction of the map \eqref{PhiCn}
to any ball $B_\rho = \{ z \in \C^n \ | \ \| z \| < \rho \}$
is convex, and it is open as a map to its image.

\begin{proof}
Consider the following commuting diagram of continuous maps:
\begin{equation} \labell{diagramCn}
 \xymatrix{%
  \R_+^n \times (S^1)^n \ar[rrrd]^{\text{projection}}
     \ar[d]_{(s_1,\ldots,s_n,e^{i\theta_1},\ldots,e^{i\theta_n}) 
          \mapsto 
           (s_1^{1/2}e^{i\theta_1},\ldots,s_n^{1/2}e^{i\theta_n})} &&& \\
 \C^n \ar[rrr]_{(z_1,\ldots,z_n) \mapsto (|z_1|^2,\ldots,|z_n|^2)} 
                                    &&& \R_+^n .}
\end{equation}
Because the projection map is convex and open 
and the map on the left is onto,
the bottom map is convex and open.

Because the ball $B_\rho$
is the pre-image of a convex set (namely, it is the preimage of the set
$\{ (s_1,\ldots,s_n) \ | \ s_1 + \ldots + s_n < \rho^2 \}$),
the restriction of the map \eqref{PhiCn} to $B_\rho$ is also convex
and open as a map to its image.
\end{proof}
\end{noTitle}

\begin{noTitle}
Let $\alpha_1,\ldots,\alpha_n$ be any vectors.  Then the map 
\begin{equation} \labell{PhiH}
\Phi_H \colon (z_1,\ldots,z_n) 
  \mapsto \sum_{j=1}^n |z_j|^2 \alpha_j
\end{equation}
is convex, and it is open as a map to its image.

Moreover, the restriction of $\Phi_H$ to any ball 
$B_\rho = \{ z \in \C^n \ | \ \| z \| < \rho \}$
is convex, and it is open as a map to its image.

\begin{proof}
Because the map~\ref{PhiCn} is convex, so is its composition 
with the linear map $(s_1,\ldots,s_n) \mapsto 
 (s_1 \alpha_1 + \ldots + s_n \alpha_n)$.
Because the restriction of a linear map to the positive orthant $\R_+^n$ 
is open as a map to its image\footnote{
This is a consequence of the following lemma:

\begin{quotation}
For any vectors $\alpha_1,\ldots,\alpha_n \in \R^{\ell}$
there exists $\eps > 0$ such that for every $\beta = \sum s_j \alpha_j$
with all $s_j \geq 0$, if $\| \beta \| < \eps$
then there exists $s' = (s'_1,\ldots,s'_n)$ such that
$\beta = \sum s'_j \alpha_j$ and $\| s' \| < 1$.
\end{quotation}

\begin{proof}[Proof of the lemma:]
Let $\beta = \sum s_j \alpha_j$ with all $s_j \geq 0$.
Then there exist $s'_j$ such that $\beta = \sum s'_j \alpha_j$,
all $s'_j \geq 0$, and the vectors $\{ \alpha_j \ | \ s'_j \neq 0 \}$
are linearly independent; cf.\ Carath\'eodory's theorem 
in convex geometry.
Let $J = \{ j \ | \ s'_j \neq 0 \}$.
The map $s \mapsto \sum s_j \alpha_j$ from $\R^J$ to
$\text{span} \{ \alpha_j \ | \ j \in J \}$
is a linear isomorphism; denote its inverse by $L_J$.
Then $s' = L_J(\beta)$, so
$\| s' \| \leq \| L_J \| \| \beta \|$
where $\| L_J \|$ is the operator norm.
The lemma holds with any $\eps < \min\limits_J
\left\{ \frac{1}{\| L_J \|} \right\}$
where $J$ runs over the subsets of $\{ 1,\ldots,n \}$
for which $\{ \alpha_j \ | \ j \in J \}$
are linearly independent.
\end{proof}
}, 
so is this composition. 
Because the map~\eqref{PhiH} is open as a map to its image,
so is its restriction to the open ball $B_\rho$.
Because this restriction is the composition of a convex map 
with a linear projection, it is convex.
\end{proof}
\end{noTitle}

We proceed with applications to symplectic geometry.
Relevant definitions can be found e.g.\ in the original paper \cite{GS}
of Guillemin and Sternberg.
We first describe local models for Hamiltonian torus actions.

\begin{noTitle} \labell{modelY}
Let $T \cong (S^1)^k$ be a torus, $\t \cong \R^k$ 
its Lie algebra, and $\t^* \cong \R^k$ the dual space.
Let $H \subset T$ be a closed subgroup, $\h \subset \t$ its Lie algebra,
and $\h^0 \subset \t^*$ the annihilator of~$\h$ in~$\t^*$.
Let $H$ act on $\C^n$ through a group homomorphism $H \to (S^1)^n$
followed by coordinatewise multiplication.
The corresponding quadratic moment map, $\Phi_H \colon \C^n \to \h^*$,
has the form $z \mapsto \sum_{j=1}^n |z_j|^2 \alpha_j$
where $\alpha_1,\ldots,\alpha_n$ are elements of $\h^*$
(namely, they are the weights of the $H$ action on $\C^n$,
times $\frac{1}{2}$).

Consider the model 
$$  Y = T \times_H \C^n \times \h^0;$$
its elements are represented by triples $[a,z,\nu]$
with $a \in T, z \in \C^n$, and $\nu \in \h^0$,
with $[ab,z,\nu] = [a, b\cdot z, \nu]$ for all $b \in H$.
Fix a splitting $\t^* = \h^* \oplus \h^0$, and consider the map
$$\Phi_Y \colon T \times_H \C^n \times \h^0 \to \t^* 
\quad , \quad
 [a,z,\nu] \mapsto \Phi_H(z) + \nu .$$

The map $\Phi_Y$ is convex and is open as a map to its image.
This follows from the commuting diagram
$$ \begin{CD}
T \times \C^n \times \h^0 @> (a,z,\nu) \mapsto (\Phi_H(z),\nu) >> 
         \h^* \times \h^0 \\
@VVV @V \cong VV \\
T \times_H \C^n \times \h^0 @> \Phi_Y >> \t^*,
\end{CD}$$
in which the top map is convex and is open as a map to its image,
the map on the left is onto, and the map on the right
is a linear isomorphism.

Similarly, if $D \subset \C^n$ and $D' \subset \h^0$
are disks centered at the origin, the restriction of $\Phi_Y$
to the subset $T \times_H D \times D'$ of $T \times_H \C^n \times \h^0$
is convex and is open as a map to its image.
This follows from the diagram
$$\begin{CD}
T \times D \times D' @> (a,z,\nu) \mapsto (\Phi_H(z),\nu) >>
            \h^* \times \h^0 \\
@VVV @V \cong VV \\
 T \times_H D \times D' @> \Phi_Y >> \t^*.
\end{CD}$$
\end{noTitle}

\begin{Proposition} \labell{GS local convexity}
Let $T$ act on a symplectic manifold with a moment map 
$\Phi \colon M \to \t^*$.
Then each point of $M$ is contained in an open set $U \subset M$
such that the restriction of $\Phi$ to~$U$ is convex
and is open as a map to its image, $\Phi(U)$.
\end{Proposition}

\begin{proof}
Fix a point $x \in M$.  

There exists a $T$-invariant neighbourhood $U$ of $x$ 
and an equivariant diffeomorphism
$f \colon U \to T \times_H D \times D'$
that carries $\Phi|_U$ to a map that differs from $\Phi_Y$
by a constant in $\t^*$,
where the model $T \times_H D \times D'$ and the map $\Phi_Y$
are as in \S\ref{modelY}.
This follows from the local normal form theorem
for Hamiltonian torus actions \cite{GS2}.
Because the restriction of $\Phi_Y$ to $T \times_H D \times D'$ 
is convex and is open as a map to its image,
so is $\Phi|_U$.
\end{proof}

We can now recover the convexity theorem of Atiyah,
Guillemin, and Sternberg along the lines given by
Condevaux-Dazord-Molino.

\begin{Theorem}
Let $M$ be a manifold equipped with a symplectic form and a torus action,
and let $\Phi \colon M \to \t^*$ be a corresponding moment map.
Suppose that $\Phi$ is proper as a map to some convex subset of $\t^*$.
Then the image of $\Phi$ is convex, its level sets are connected,
and the moment map is open as a map to its image.
\end{Theorem}

\begin{proof}
By Proposition \ref{GS local convexity},
every point in $M$ is contained in an open set $U$
such that the map $\Phi|_U$ is convex and is open as a map
to its image, $\Phi(U)$.  The conclusion then follows 
from Theorem \ref{Theorem}. 
\end{proof}

\section{The results of Birtea-Ortega-Ratiu}
\labell{sec:ratiu}

The paper \cite{BOR1} of Birtea, Ortega, and Ratiu
contains results that are similar to ours.
For the benefit of the reader, we present their results here. 

\begin{Theorem}[{\cite[Theorem 2.28]{BOR1}}]
Let $X$ be a topological space that is connected, locally connected, 
first countable, and normal. Let $V$ be a finite dimensional vector space.
Let $f \colon X \to V$ be a map that satisfies the following conditions.
\begin{enumerate}
\item
The map $f$ is continuous and is closed.
\item
The map $f$ has \emph{local convexity data}:
for each $x \in X$ and each sufficiently small neighborhood $U$
of $x$ there exists a convex cone $C \subset V$ with vertex
at $f(x)$ such that the restriction $f|_U \colon U \to C$
is an open map with respect to the subset topology on $C \subset V$.
\item
The map $f$ is \emph{locally fiber connected}:
for each $x \in X$, any open neighborhood of $x$
contains a neighborhood $U$ of $x$ that does not
intersect two connected components of the fiber
$f\inv(f(x'))$ for any $x' \in U$.
\end{enumerate}
Then the fibers of $f$ are connected, the map $f$ is open onto its image,
and the image $f(X)$ is a closed convex set.
\end{Theorem}

\begin{Remark}
The paper \cite{BOR2} contains a more general convexity result;
in particular it contains a more liberal definition of 
having local convexity data:  for each $x \in X$ there exist
arbitrarily small neighborhoods $U$ of $x$ such that $f(U)$ is convex
\cite[Def.2.8]{BOR2}.  Here, openness of the maps $f|_U \colon U \to f(U)$
is not part of the definition of ``local convexity data", 
but it is assumed separately.
\end{Remark}

Birtea-Ortega-Ratiu also sketch a proof of the following
infinite dimensional version:

\begin{Theorem}[{\cite[Theorem 2.31]{BOR1}}]
Let $X$ be a topological space that is connected, locally connected, 
first countable, and normal. Let $(V,\| \ \|)$ be a Banach space 
that is the dual of another Banach space.  Let $f \colon X \to V$ 
be a map that satisfies the following conditions.
\begin{enumerate}
\item
The map $f$ is continuous with respect to the norm topology on $V$ 
and is closed with respect to the weak-star topology on $V$.  
\item
The map $f$ has \emph{local convexity data} (see above).
\item
The map $f$ is \emph{locally fiber connected} (see above).
\end{enumerate}
Then the fibers of $f$ are connected,
the map $f$ is open onto its image with respect to the weak-star
topology, and the image $f(X) \subset V$ is convex and is closed 
in the weak-star topology.
\end{Theorem}


\begin{Remark}\ 

\begin{itemize}

\item
We work with a convex subset of $V$; they similarly note 
that their theorem remains true 
with $V$ replaced by a convex subset of $V$
\cite[remark 2.29]{BOR1}.

\item
In \cite{BOR2} they allow more general target spaces,
which are not vector spaces.

\item
We assume that the domain is Hausdorff and the map is proper
(in the sense that the preimage of a compact set is compact);
they assume that the domain is first countable and normal
and that the map is closed.  
We are not aware of non-artificial examples where one of these
assumptions holds and the other doesn't.

\item
We assume that each point is contained in an open set
on which the map is a \emph{convex map},
a condition that we define in Definition~\ref{def:convex}, 
section~\ref{sec:local-global}.
They assume that the map \emph{has local convexity data}
(defined in \cite[Def.2.7]{BOR1} and re-defined in \cite[Def.2.8]{BOR2})
and satisfies the \emph{locally fiber connected condition}
(defined in \cite[Def.2.15]{BOR1} as a slight generalization
of \cite[\S3.4, after Def.3.6]{benoist}).

\begin{Example*}
The inclusion map of a closed ball into $\R^n$ 
is a convex map in our sense.
It does not have local convexity data in the sense of \cite{BOR1},
but it does have local convexity data in the sense of \cite{BOR2}.
\end{Example*}

\item
If a map is convex, then it has local convexity data
(in the broader sense, of \cite{BOR2})
and it is locally fiber connected.
Thus, our ``convexity/connectedness" assumptions are stricter
than those of \cite{BOR1}, but our conclusion is stronger.

\item
Both we and \cite{BOR1} allow a broad interpretation of ``local":
\begin{itemize}
\item
In \cite{BOR1}, the ``locally fiber connected condition"
on a subset $A$ of $X$ with respect to a map $f \colon X \to V$
depends not only on the restriction of the map $f$ to the set $A$
but also on the information of which points in $A$
belong to the same fiber \emph{in $X$}.
(This is where the definition of \cite{BOR1} differs
from that of Benoist.)
\item
In our paper, we assume that each point is contained in an open set
on which the map is convex and is open as a map to its image,
but we do not insist that these open sets form a basis
to the topology (cf.\ Remark~\ref{U not small}).
(E.G., in the presence of a group action, it's fine to just check
neighborhoods of orbits rather than neighborhoods of individual
points.)
\end{itemize}

\end{itemize}
\end{Remark}

\end{document}